\documentclass{elsarticle}
\usepackage{amsmath,amssymb,amsthm}
\usepackage{algorithm}
\usepackage[noend]{algpseudocode}

\usepackage{ccicons}
\usepackage{hyperref}

\usepackage[show]{ed}

\makeatletter
\newcommand{\bigger}{\bBigg@{3}}
\newcommand{\vast}{\bBigg@{4}}
\newcommand{\Vast}{\bBigg@{5}}
\makeatother

\newcommand{\rd}{\mathrm{d}}
\newcommand{\calD}{\mathcal{D}}

\newcommand{\calO}{\mathcal{O}}

\newcommand{\calZ}{\mathcal{Z}}

\newcommand{\bbE}{\mathbb{E}}
\newcommand{\bbN}{\mathbb{N}}
\newcommand{\bbP}{\mathbb{P}} 
\newcommand{\bbZ}{\mathbb{Z}}
\newcommand{\rth}{\mathrm{th}}

\makeatletter
\def\ps@pprintTitle{%
 \let\@oddhead\@empty
 \let\@evenhead\@empty
 \def\@oddfoot{\footnotesize\textcopyright 2017. This manuscript licensed under \href{http://creativecommons.org/licenses/by-nc-nd/4.0/}{Creative Commons CC-BY-NC-ND 4.0} \ccbyncnd}%
 \let\@evenfoot\@oddfoot}
\makeatother

\begin{document}
\begin{frontmatter}
\title{Subdiffusive discrete time random walks via Monte Carlo and subordination}
\author[SU]{J.A. Nichols}
\ead{nichols@ljll.math.upmc.fr}
\author[UNSW]{B. I. Henry}
\ead{b.henry@unsw.edu.au}
\author[UNSW]{C. N. Angstmann\corref{cor}}
\ead{c.angstmann@unsw.edu.au}
\cortext[cor]{Corresponding author}
\address[SU]{Sorbonne Universit\'{e}s, UPMC Univ Paris 06, CNRS, UMR 7598, Laboratoire Jacques Louis-Lions, 4, place Jussieu 75005, Paris France}
\address[UNSW]{School of Mathematics and Statistics, University of New South Wales, Sydney NSW 2052 Australia}

\date{\today}
\begin{abstract}
A class of discrete time random walks has recently been introduced to provide a stochastic process based numerical scheme for solving fractional order partial differential equations, including the fractional subdiffusion equation.
Here we develop a Monte Carlo method for simulating discrete time random walks with Sibuya power law waiting times, providing
another approximate solution of the fractional subdiffusion equation. 
The computation time scales as a power law in the number of time steps with a fractional exponent simply related to the order of the fractional derivative.
We also provide an explicit form of a subordinator for discrete time random walks with Sibuya power law waiting times. This subordinator transforms from an operational time, in the expected number of random walk steps, to the physical time, in the number
of time steps.
\end{abstract}



\end{frontmatter}

\section{Introduction}

There has been a great deal of interest in recent years in solving partial differential equations (PDEs) with fractional order derivatives. This interest has been motivated by the recognition that physical properties, such as anomalous 
diffusion \cite{MK2000} and  viscoelasticity \cite{M2010}, can be modelled using fractional order PDEs. A classic example
is the time fractional subdiffusion equation, 
\begin{equation} \label{eq:frac_diffusion_intro}
 \frac{\partial u}{\partial t} = D_\alpha \,_0\mathcal{D}_t^{1-\alpha} \frac{\partial^2 u}{\partial x^2},
\end{equation}
where $_0\mathcal{D}_t^{1-\alpha}$ denotes the Riemann-Liouville fractional derivative of order $1-\alpha$.
This fractional PDE governs the evolution of the probability density for an ensemble of diffusing particles
whose variance scales as a 
sublinear power law with time,
\begin{equation}
\bbE[\| x(t) - x(0) \|^2] \sim t^\alpha, \, \quad 0<\alpha<1.
\end{equation}
Brownian motion, described by the standard diffusion equation, is recovered when $\alpha=1$.
The time fractional subdiffusion equation has been derived as a limit form of the evolution equation for the probability density
describing an ensemble of continuous time random walks (CTRWs), characterized by power law waiting time densities and nearest neighbour, or Gaussian,
step length densities \cite{MK2000}. The power law waiting time density models the phenomenon of trapping whereby the longer a particle remains at a site the more likely it is to continue waiting.
The derivation, using CTRWs with power law waiting times, has also been extended to fractional Fokker-Plank equations \cite{BMK2000,SK2006,HLS2010} and fractional reaction-diffusion equations \cite{HLW2006,F2010}.

Algebraic approaches to solving fractional PDEs include the variational iteration method \cite{H1998}, the
homotopy perturbation method \cite{MO2007}, and the Admonian decomposition method \cite{DRBW2012}. Numerical approaches include explicit \cite{YA2005} and implicit \cite {LH2005} finite difference schemes, Galerkin methods \cite{MM2009}, and spectral methods \cite{CLTA2007}.
In earlier work, inspired by the CTRW method of deriving fractional PDEs, we introduced a numerical method for approximating solutions of fractional PDEs based on discrete time random walks (DTRWs) \cite{ADHN2015, ADHJLN2016}. The DTRWs are characterized by a discrete waiting time probability mass function and a step-length probability mass function. With appropriately chosen, waiting time and step length, mass functions the generalized master equation for the evolution
of the probability mass function, describing the position and location of an ensemble of DTRWs, limits to the same generalized master equation for the ensemble of CTRWs with power law waiting time densities. This correspondence can be exploited to provide an explicit finite difference scheme whose solution approximates the solution of the fractional PDE. This DTRW method can also be used to approximate solutions to fractional order ODE compartment models \cite{AHM2016}, and integer order advection-diffusion PDEs \cite{AHJM2016}.

In this work we develop a Monte Carlo method for simulating DTRWs \cite{ADHN2015,ADHJLN2016} that provides a numerical method 
for approximating the solution of the fractional subdiffusion equation. 
Monte Carlo methods based on   
simple discrete random walks have been widely employed in solving diffusion problems
in a variety of application areas, from mathematical finance \cite{Boyle1977, KS1998}, to the simulation of chemical reactions \cite{Gillespie1977}, to mathematical biology \cite{CPB2008}. 

In the DTRW that we introduced for anomalous subdiffusion \cite{ADHN2015,ADHJLN2016} the waiting time between jumps is a simple renewal process \cite{Doob1948}, and the waiting time probability mass function is
 the Sibuya distribution \cite{Sibuya79, Sibuya80} given by 
\[
 \bbP[\text{Jump waiting time} = m ] = \frac{\alpha}{m} \prod_{\ell=1}^{m-1} \left(1 - \frac{\alpha}{\ell}\right) .
\]
The combination of this waiting time probability mass function together with nearest-neighbour jumps yields a process that is subdiffusive, in the diffusion limit of the space-time grid. 

Monte Carlo path-wise simulations for calculating the ensemble average solution of the subdiffusive process out performs finite difference methods for small $\alpha$, typically less than $3/4$. The performance benefit is due to the semi-Markovian nature of the process. In the Monte Carlo simulations we do not need to know the history of a particle when choosing its next jump. Furthermore we find that the expected number of jumps for a trajectory up to a time point $n$ scales as $\calO(n^\alpha)$, meaning less jumps are required for small $\alpha$. We find a closed form expression for the probability of a process undergoing $k$ jumps up to time $n$, 
 \begin{equation*}
 \bbP[ \text{Exactly $k$ jumps up to time $n$} ] = (-1)^n 
\sum_{\ell = 0}^{k} (-1)^\ell \binom{k}{\ell}  \binom{(\ell+1)\alpha - 1}{n} 
\, .
\end{equation*}
We use this expression to calculate expected computational complexity for a path-wise Monte Carlo simulation of our process, and compare it with the corresponding finite difference scheme that calculates the ensemble average result directly. 
The expression above furthermore enables us to find a closed form solution for the ensemble average, in a formula that is reminiscent of the subordination formula used to solve the subdiffusive time-fractional partial differential equation \cite{B2001}.

The remainder of this paper is as follows. In Section 2 we derive the generalized master equation for  the ensemble average of discrete time random walk paths characterized by
independently distributed waiting times and jump lengths.
The particular cases of exponentially distributed and Sibuya distributed waiting times are discussed, leading to the standard diffusion equation, and the fractional diffusion equation, respectively. 
In Section 3 we describe the Monte Carlo method for approximating the solution of the fractional diffusion equation.
In Section 4 we compute the computational complexity of calculating  an ensemble of random paths.
In Section 5 we derive the discrete subordination result linking the solution of the generalized master equation, with Sibuya distributed waiting times,
to the solution of the generalized master equation with exponentially distributed waiting times.
Finally, in Section 6 we present numerical results of the Monte Carlo path-wise simulations and make comparisons to known analytic solutions, as well as to numerical results based on the corresponding finite difference scheme. 
We show that CPU timing of the Monte Carlo simulations is consistent with the computational complexity analysis in Section 4.

\section{Generalized master equations for an ensemble of discrete time random walks}

We consider a random walker, which we refer to as a particle, to be traversing discrete space time. The particle undergoes a series of jumps from point to point. After each jump we draw a random waiting time, and after waiting at that point for that amount of time, the particle instantaneously jumps to some other point, chosen randomly from a jump probability distribution. Thus we generate a path by producing a series of independently distributed waiting times $W_k(\omega)\in \bbN$ and jumps $J_k(\omega)\in\bbZ$. The process can easily be  generalised to $d$-dimensional space by allowing $J_k(\omega)\in\bbZ^d$, however for simplicity we only consider one spatial dimension here. The waiting time, and jumps, are both discrete renewal processes with independent probability mass functions, their combination forms a discrete time random walk. The index $k\in\bbN$ denotes the jump index, $J_k$ denotes the length and direction of the $k^\rth$ jump, and $W_k(\omega)$ denotes the time the particle waits between the $(k-1)^\rth$ and $k^\rth$ jump. The variable $\omega\in\Omega$ is our realisation for the random variables with distribution $\gamma$. 

For simplicity we assume the particle starts at a point $i_0$ and at time $0$. Thus after $k$ jumps, the particle will be at a position 
\begin{equation} \label{eq:def_X_k}
 X_k(\omega) = i_0 + \sum_{\ell=1}^k J_\ell(\omega) \, ,
\end{equation}
with $X_0(\omega) = i_0$, and the $k^\rth$ jump will occur at the time
\begin{equation} \label{eq:def_T_k}
 T_k(\omega) = \sum_{\ell=1}^k W_\ell(\omega) \, ,
\end{equation}
with $T_0(\omega) = 0$.

To find the position of the particle at a particular time $n\in\bbN$, we first count the number of jumps up to that time, defining the count as 
\[
k_n(\omega) = \max\{ k \,:\, T_k(\omega) \le n \}. 
\]
Whatever time $n$ we choose, we know that  $T_{k_n(\omega)}(\omega) \le n < T_{k_{n+1}(\omega)}(\omega)$. Thus we can define the position of the particle in terms of time $n$ as 
\[
 X(n,\omega) := X_{k_n(\omega)}(\omega) = i_0 + \sum_{k=1}^{k_n(\omega)} J_k(\omega)
\]

The motion of the particle is determined by the distribution of the jumps and waiting times. The probability of jumping a distance $j$, conditional on the particle being at position $i$ when the jump occurs at time $n$ is denoted by
\[
\lambda(j;\, i,n) := \bbP_\gamma[J_k(\cdot) = j \text{ for any $k$}\,|\, X_{k-1}(\cdot) = i \text{ and } T_{k}(\cdot) = n] \, . 
\]
 We also specify the waiting time probability mass function
 \begin{equation} \label{eq:def_phi}
\phi(m) := \bbP_\gamma[W_k(\cdot) = m \text{ for any $k$}].
\end{equation}
Note that only a single jump is allowed in any one time step and so $\phi(0)=0$. This could readily be extended to space and time dependent waiting time probability mass functions $\phi(m;\, i, n)$ but we have not pursued this further here. 
The probability mass functions satisfy the normalization conditions $\sum_j \lambda(j;\,i,n) = 1$ and $\sum_{m\ge 0} \phi(m) = 1$. It is the specific choice of the mass functions that yields different types of diffusive movement. We also make use of the \emph{survival probability}, defined as the likelihood of no jump occurring up until $m$ time steps, 
\begin{equation} \label{eq:def_survival}
\Phi(m) := \bbP_\gamma[W_k(\cdot) > m \text{ for any $k$} ] = 1 - \sum_{\ell=1}^m \phi(\ell) \, ,
\end{equation}
Evidently $\Phi(m)\to 0$ as $m\to\infty$. Furthermore $\Phi(m)$ is a decreasing function as all $\phi(\ell)$ are non-negative.

In the remainder of this article we consider nearest neighbour jumps, i.e.
\begin{equation} \label{eq:nn_jump_prob}
 \lambda(j-i;\, i, n) =  p_\ell \delta_{j, i+1} + p_r \delta_{j, i-1}  \, ,
\end{equation}
where $p_r, p_\ell \in [0,1]$ are the right and left jump probabilities respectively, and $p_r = 1-p_\ell$.
We also consider a slight generalization, allowing for self jumps to occur, with some likelihood $r\in(0,1]$, viz;
\begin{equation} \label{eq:nn_r_jump_prob}
 \lambda(j-i;\, i, n) = r \left( p_\ell \delta_{j, i+1} + p_r \delta_{j, i-1} \right) + (1-r) \delta_{i,j} \, .
\end{equation}
External forces could be simulated by considering  biased nearest neighbour jumps but we have not considered this
further here. 

\subsection{The ensemble average solution}

We can obtain a ensemble average solution of this stochastic process, which gives us the distribution of the location of the particle through space and time. That is, we wish to find the distribution 
\begin{equation} \label{eq:U_def}
U(i,n) = \bbP_\gamma [ X(n,\cdot) = i ].
\end{equation}

 In the tradition of continuous time random walks, we derive two quantities, the arrival density, which gives us the likelihood that the particle will jump to a particular point at a particular time, and the location density, which is the likelihood of a particle being at an exact point at a certain time. These two quantities allow us to keep track of occupancy times between jumps, hence allowing us to consider anomalous diffusion, where the process is Markovian in terms of jumps, but non-Markovian in time.

The arrival density $Q_k(i,n)$ is defined as the likelihood of the particle arriving at location $i$ at exactly time $n$, after having completed $k$ jumps. Formally we can write,
\begin{equation} \label{eq:def_Q_k}
 Q_k(i,n) := \bbP_\gamma \left [ X_k(\cdot) = i \text{ and } T_k(\cdot) = n \right ] \, .
\end{equation}
The arrival density at the $k^\rth$ jump can be calculated from the arrival density for the previous jump, we expand \eqref{eq:def_Q_k} and using Bayes' law obtain 
\begin{align}
 Q_k(i,n) 
&= \sum_{j} \bbP_\gamma  [  X_{k-1}(\cdot) = j \text{ and } J_k(\cdot) = i-j \text{ and } T_{k}(\cdot) = n ] \nonumber \\
&= \sum_{j} \bbP_\gamma  [ J_k(\cdot) = i-j \,|\, X_{k-1}(\cdot) = j \text{ and } T_{k}(\cdot) = n ] \times \bbP_\gamma[ X_{k-1}(\cdot) = j \text{ and } T_{k}(\cdot) = n ] \nonumber \\
&= \sum_{j} \bbP_\gamma  [ J_k(\cdot) = i-j \,|\, X_{k-1}(\cdot) = j \text{ and } T_{k}(\cdot) = n ] \nonumber \\ 
& \quad \times\sum_{m=0}^n \bbP_\gamma[ X_{k-1}(\cdot) = j \text{ and } T_{k-1}(\cdot) = m ] \times \bbP_\gamma[  W_{k}(\cdot) = n-m ] \nonumber \\
&= \sum_{j} \sum_{m=0}^n  \lambda(i-j;\, j, n) \, \phi(n - m) \, Q_{k-1}(j, m) \, . \label{eq:recursive_Qk}
\end{align}
The arrival density after any arbitrary number of jumps is simply
\[
 Q(i,n) := \sum_{k=0}^\infty Q_k(i,n) \, ,
\]
and summing both sides of \eqref{eq:recursive_Qk} for $k=1,\ldots,\infty$, we obtain
\begin{align}
 Q(i,n) 
 &= Q_0(i,n) + \sum_{j} \sum_{m=0}^n  \lambda(i-j;\, j, n) \, \phi(n - m) \, Q(j, m) \nonumber \\
 &= \delta_{i,i_0}\delta_{n,0} + \sum_{j} \sum_{m=0}^n  \lambda(i-j;\, j, n) \, \phi(n - m) \, Q(j, m) \, , \label{eq:Q_recursion}
\end{align}
where $Q_0(i,n)$ is the very first arrival, or the initial condition, which we see is equal to $\delta_{i,i_0} \delta_{n,0}$ from \eqref{eq:def_X_k}, \eqref{eq:def_T_k} and \eqref{eq:def_Q_k}.

The distribution $U(i,n)$ can now be obtained from the arrival probabilities, viz;
\begin{align}
 U(i,n) 
&= \bbP_\gamma[X(n,\cdot) = i] \nonumber \\
&=\sum_{k=0}^\infty \bbP_\gamma[ X_k(\cdot) = i \text{ and } T_k(\cdot) \le n \text{ and } W_{k+1}(\cdot) > n - T_k(\cdot) ] \nonumber \\
&=\sum_{m=0}^{n} \sum_{k=0}^\infty \bbP_\gamma[ X_k(\cdot) = i \text{ and } T_k(\cdot) = m ] \times \bbP_\gamma[W_{k+1}(\cdot) > n - m ] \nonumber \\
&=\sum_{m=0}^{n} Q(i,m) \, \Phi(n - m) \label{eq:U_Q_dep} 
\end{align}

To derive a \emph{generalised master equation} we calculate the change in $U$ between adjacent time-steps. Using \eqref{eq:Q_recursion} and \eqref{eq:U_Q_dep} and noting that $\Phi(n-m) = \Phi(n-1-m) - \phi(n-m)$, we obtain
\begin{align} 
U(i,n) - U(i,n-1) 
&= Q(i,n) + \sum_{m=0}^{n-1} Q(i,m) \, (\Phi(n - m) - \Phi(n-1-m)) \nonumber \\
&= Q(i,n) - \sum_{m=0}^{n-1} Q(i,m) \, \phi(n - m) \nonumber \\
&= \delta_{i,i_0}\delta_{n,0} + \sum_{j} \sum_{m=0}^n  \lambda(i-j;\, j, n) \, \phi(n - m) \, Q(j, m) - \sum_{m=0}^{n-1} Q(i,m) \, \phi(n - m) \nonumber \\
&= \delta_{i,i_0}\delta_{n,0} + \sum_{j}(\lambda(i-j;\, j, n) - \delta_{i,j})  \sum_{m=0}^n \phi(n - m) \, Q(j, m). \label{eq:GME_Q}
\end{align}
We now make use of the $\calZ$-transform, see for example \cite{OSB1989},
\[
 \calZ_n\{ F ;\, z\} = \sum_{n=0}^\infty F(n) \, z^n \, .
\]
To remove the dependence on the arrival density $Q$ in \eqref{eq:GME_Q}, we use the convolution property of the $\calZ$-transform. Applying the $\calZ$-transform to \eqref{eq:U_Q_dep} we obtain
\begin{equation} \nonumber
\calZ_n\{ U(i,n) ;\, z\} = \calZ_n\left\{ \sum_{m=0}^n \Phi(n - m) \, Q(i,m) ;\, z\right\} = \calZ_n\{ Q(i,n) ;\, z\} \calZ_n\{ \Phi(i,n) ;\, z\}.
\end{equation}
We now use this result in the convolution in \eqref{eq:GME_Q} to obtain
\begin{align} \label{eq:U_Q_rel}
 \calZ_n\left\{ \sum_{m=0}^n \phi(n - m) \, Q(i,m);\, z\right\} 
 &= \calZ_n\{ Q(i,n) ;\, z\} \calZ_n\{ \phi(i,n) ;\, z\}  \nonumber\\ 
 &= \calZ_n\{ U(i,n) ;\, z\} \frac{\calZ_n\{ \phi(i,n) ;\, z\} }{ \calZ_n\{ \Phi(i,n) ;\, z\}}  \nonumber \\
 &= \calZ_n\{ U(i,n) ;\, z\} \calZ_n\{ K(n) ;\, z \} ,
\end{align}
where we have have introduced $K$, which we call the \emph{memory kernel}. It is calculated from the waiting time probability mass function as
\begin{equation}
 \calZ_n\{ K(n) ;\, z \} = \frac{\calZ_n\{\phi(n) ;\, z\} }{\calZ_n\{\Phi(n) ;\, z\}} \, .
\end{equation}
We can then use  the convolution property of the $\calZ$-transform again to rewrite \eqref{eq:GME_Q} in terms of $K$ for $n\geq1$, 
\begin{equation} \label{eq:GME}
 U(i,n) - U(i,n-1) = \sum_j \left( \lambda(i-j;\, j, n) - \delta_{i,j} \right) \sum_{m=0}^n K(n - m) \, U(j, m) \, .
\end{equation}

Equation \eqref{eq:GME} describes a finite difference scheme that propagates the probability density for the stochastic process.  The master equation in \eqref{eq:GME}, and other related master equations, were derived previously \cite{ADHN2015,ADHJLN2016}  using a slightly less formal approach.
\subsection{Exponentially distributed waiting times}

First we consider the simple case of a Markovian process, where the likelihood of a jump occurring is independent of when the particle arrived at its current location. For further simplicity we consider a fixed probability $\sigma \in[0,1]$ of a jump occurring in a time step, so the probability of not jumping for $m$ steps is $(1-\sigma)^m$, hence our survival probability and waiting time probability mass functions are
\[
 \Phi(m) = (1-\sigma)^m \quad\text{and}\quad \phi(m) = \sigma(1-\sigma)^{m-1}, 
\]
respectively. In this case the kernel is given by
\[
K(m) = \sigma \delta_{1, m}.
\]
Using this form of the kernel in \eqref{eq:GME}  with nearest-neighbour jumps \eqref{eq:nn_r_jump_prob}, and with equal probability of left or right jumps (that is $p_r = p_\ell = 1/2$), we arrive at the discrete stencil for obtaining an approximate solution to the classic diffusion equation, 
\begin{equation}
 U(i,n) - U(i,n-1) = \frac{\sigma r}{2}\big[ U(i-1, n-1) - 2U(i,n-1) + U(i+1, n-1)\big] \, .
\end{equation}

\subsection{Power law distributed waiting times} \label{sec:power_law}

We now  consider a discrete power-law waiting time distribution, resulting in a discrete process that in the continuum limit approaches anomalous diffusion \cite{ADHN2015,ADHJLN2016}. We have a parameter $0\le \alpha \le 1$ that gives us the exponent of the anomalous diffusion that we require. We take a process given by the jump probability in \eqref{eq:nn_jump_prob}, and the Sibuya distribution \cite{Sibuya79} for the waiting time probability mass function
\begin{equation} \label{eq:sibuya_def}
 \phi(m) = \frac{\alpha}{m} \prod_{\ell=1}^{m-1} \left( 1 - \frac{\alpha}{\ell} \right)  \quad\text{with}\quad \phi(0) = 0.
\end{equation}
 We can see quite clearly that this waiting time probability mass function can be decomposed into two components, the probability that the jump will occur at exactly time $m$, given by $\alpha / m$, and the probability that no jump will occur up until  time $m$, given by $\prod_{\ell=1}^{m-1} \left( 1 - \alpha / \ell \right)$. In fact we have from \eqref{eq:def_survival} that 
\begin{equation} \label{eq:sibuya_surv_def}
\Phi(m) = \prod_{\ell=1}^{m} \left( 1 - \frac{\alpha}{\ell} \right) = \binom{m-\alpha }{m} \, .
\end{equation}
We note that the waiting time probability mass function can also be expressed in terms of binomial coefficients,
\[
 \phi(m) =-\binom{m-1-\alpha }{m} - \delta_{0,m} \, .
\]
These forms enable easy calculation of the $\calZ$-transforms and hence the memory kernel, explicitly found as 
\begin{equation} \label{eq:sibuya_kernel}
 K(m) = \prod_{\ell=1}^{m} \left ( 1 - \frac{2-\alpha}{\ell} \right ) \,\text{, with } K(0)=0 \text{ and } K(1) = \alpha \, .
\end{equation}
The memory kernel also has a binomial representation,
\[
 K(m) = \binom{1-\alpha}{m} (-1)^{m} - \delta_{0,m} + \delta_{1,m} \, .
\]

When we use this kernel in the generalised master equation \eqref{eq:GME}, with nearest-neighbour jumps and self-jumps, as specified in \eqref{eq:nn_r_jump_prob}, we obtain
\begin{equation} \label{eq:sibuya_GME}
\begin{split}
 U(i,n) - U(i,n-1) = \sum_{m=0}^{n-1} \left( \binom{1-\alpha}{n-m} (-1)^{n-m} + \delta_{1,n-m} \right) \\
  \quad \times \frac{r}{2}\big(U(i-1, m) - 2U(i,m) + U(i+1, m)\big) \, .
\end{split}
\end{equation}
Equation \eqref{eq:sibuya_GME} provides a finite-difference  approximation to the fractional diffusion equation. 
Specifically, taking
\begin{equation}
\lim_{\Delta x, \Delta t \to 0} U(\lfloor x/\Delta x\rfloor, \lfloor t/\Delta t \rfloor) = u(x,t),
\end{equation}
in \eqref{eq:sibuya_GME}, using Taylor expansions for nearest neighbour terms, and taking the limit
$\lim_{\Delta x, \Delta t \to 0}$ we find that
$u(x,t)$ is the solution of the fractional diffusion equation,
\begin{equation} \label{eq:frac_diffusion}
 \frac{\partial u}{\partial t} = D_\alpha \,_0\calD_t^{1-\alpha} \frac{\partial^2 u}{\partial x^2},
\end{equation}
where $\,_0\calD_t^{1-\alpha}$ is the Riemann-Liouville fractional derivative, and where we have defined the diffusion coefficient,
\begin{equation} \label{eq:D_alpha}
D_\alpha = \lim_{\Delta x,\Delta t\to 0} \frac{r \Delta x^2}{2 \Delta t^\alpha} < \infty .
\end{equation}
The diffusion coefficient is required to remain finite in the limit $\lim_{\Delta x, \Delta t \to 0}$.

We note that the coefficients of the memory kernel in \eqref{eq:sibuya_kernel} are related to the Gr\"unwald-Letnikov coefficients for the finite-difference scheme of calculating fractional derivatives. More precisely, the backwards Gr\"unwald-Letnikov derivative of order $1-\alpha$ is given by
\[
\,_0\calD_t^{1-\alpha, \mathrm{GL}} f(t) = \lim_{\Delta t\to 0} \frac{1}{\Delta t^{1-\alpha}} \sum_{m = 0}^{\lfloor t/\Delta t \rfloor} \binom{1-\alpha}{m} (-1)^m f(t - m\Delta t)
\]

The scheme in \eqref{eq:sibuya_GME} is shown to converge to \eqref{eq:frac_diffusion} as $\Delta x, \Delta t\to 0$ in \cite{ADHJLN2016}. However $U$ can also be used as an approximation to $u$ for a fixed grid spacing. If we take some finite $\Delta x$ and $\Delta t$, then we have that $U(i,n)$ approximates $u(i\Delta x, n\Delta t)$, the solution of \eqref{eq:frac_diffusion}, where
\begin{equation} \label{eq:D_alpha_fixed}
D_\alpha = \frac{r \Delta x^2}{2 \Delta t^\alpha} \, .
\end{equation}
In practice we will want to approximate the solution to \eqref{eq:frac_diffusion} with a fixed diffusion coefficient, $D_\alpha$, and spatial grid size, $\Delta x$. This leaves us with two free parameters, $r$, which is restricted to $(0,1]$ because of its role as a probability of a jump occurring, and $\Delta t$. By fixing a value of $r$, we then have $\Delta t$ by inverting \eqref{eq:D_alpha_fixed},
\begin{equation} \label{eq:delta_t_dep}
\Delta t = \left( \frac{r \Delta x^2}{2 D_\alpha} \right)^{1/\alpha} \, .
\end{equation}

\section{Approximating the solution of the fPDE using Monte Carlo methods} \label{sec:MC}

It is straightforward to produce a path $X(n, \omega)$. To get as far as a time point $n$ we must produce $\ell = k_n(\omega)$ independent jumps $J_\ell$ and $\ell$ independent waiting times $W_\ell$  (such that $T_{k_n-1}(\omega) < n \le T_{k_n}(\omega)$). 

To produce a positive discrete random variable distributed by $\phi(m)$, we can draw a uniform random number $u\in[0,1]$ and find $m$ such that $\Phi(m) \le u < \Phi(m-1)$, recall that $\Phi$ is a decreasing function. We can calculate the correct $m$ either by knowing the inverse function $\Phi^{-1}(u)$, which plays a similar role to the inverse cummulative distribution of $\phi$, or by doing a simple search amongst a pre-computed collection of $\Phi(m)$. The jumps are nearest neighbour jumps that are straightforward to generate.

Now, suppose for some fixed time $t$ we want to approximate either the solution $u(x,t)$ of \eqref{eq:frac_diffusion} or $Au(x,t)$, where $A$ is some linear operator. If we produce $N$ individual paths $X(n, \omega_\ell)$, where $\ell=1,\ldots,N$, then we can approximate $X(i,n)$ and hence $u(x,t)$ by taking an average over the paths, viz;
\[
u(i\Delta x, n\Delta t) \sim U(i,n) \sim \frac{1}{N}  \sum_{\ell=1}^{N}[X(n, \omega_\ell) = i],
\]
where $[\cdot]$ are Iverson brackets, which are equal to $1$ when the argument is True and $0$ otherwise. Hence $\sum_{\ell=1}^{N}[X(n, \omega_\ell) = i]$ is the count of the number of paths where $X(n, \omega_\ell) = i$. 
We can now calculate expectations directly. For example, the expected position 
\begin{equation}
\bbE[x] = \int_{-\infty}^\infty x\, u(x,t) \, \rd x \sim \sum_{i=-\infty}^\infty i\, U(i,n) \sim \frac{1}{N} \sum_{\ell=1}^N  X(n, \omega_\ell).
\end{equation}

\section{Computational complexity of the Monte Carlo method for the fPDE} \label{sec:expected_jumps}
  
To know the likely computational complexity of calculating an ensemble of random paths, we need to know how many jumps are expected up to a time $n$, that is we wish to compute $\bbE_\gamma[ k : T_k \le n < T_{k+1} ]$. This is a problem encountered often in renewal process theory. For a guide to techniques such as those used here see \cite{barbu2008}.

We wish to calculate the probability of taking $k$ jumps up to and including time $n$, and taking the next jump after time $n$, that is $\bbP_\gamma[ T_{k}(\cdot) \le  n < T_{k+1}(\cdot) ]$. This is best decomposed into two quantities,
\[
   \bbP_\gamma[ T_{k}(\cdot) \le n < T_{k+1}(\cdot) ] = \bbP_\gamma[ T_{k}(\cdot) < n < T_{k+1}(\cdot) ] +  \bbP_\gamma[  T_{k}(\cdot) = n ] .
\]

For convenience we use the labels
\begin{equation} \label{eq:b_c_define}
 b(n,k) := \bbP_\gamma[ T_k(\cdot) = n ]  \quad\text{and}\quad c(n,k) := \bbP_\gamma[ T_{k-1}(\cdot) < n < T_{k}(\cdot) ] . 
\end{equation}
Using \eqref{eq:def_phi}, \eqref{eq:def_T_k} and \eqref{eq:def_survival}, it is possible to demonstrate recursion relations for both these quantities,
\begin{equation*}
b(n,k) = \sum_{m=1}^{n-1} b(m, k-1) \phi(n - m),
\end{equation*}
and note that $b(n, 1) = \bbP_\gamma[ T_1(\cdot) = n ] = \phi(n)$ for all $n$. This recursion states quite simply that the likelihood of jumping at time $n$ after $k$ jumps is the likelihood of jumping at any time $m$ before and then waiting the time $n-m$ between. We also see that we have 
\begin{equation*}
c(n, k) =
 \sum_{m=0}^{n-1}b(m, k-1)\Phi(n - m)
\end{equation*}
with $c(n, 1) = \bbP_\gamma[ T_0(\cdot) < n < T_1(\cdot) ] = \Phi(n)$ for all $n$. Note that for $n < k$ we have $b(n,k)=c(n,k+1) = 0$.

If we write $b(n, k) = \bbP_\gamma[ T_k(\cdot) = n ]$ and $c(n,k) = \bbP_\gamma[ T_{k-1}(\cdot) < n < T_{k}(\cdot) ]$ then we can see, taking the $\calZ$-transform that
\begin{align} \label{eq:bnk_induction}
\calZ_n\{b(n, k) \} 
&= \sum_{n=0}^\infty b(n, k) z^n, \nonumber \\ 
&= \sum_{n=0}^\infty \sum_{m=1}^{n-1} b(m, k-1) \phi(n - m) z^n, \nonumber \\
&= \sum_{m=1}^{\infty} b(m, k-1) z^m \sum_{n=1}^\infty \phi(n) z^{n}, \nonumber \\
&= \calZ_n\{b(n, k-1)\} \calZ_n\{ \phi(n) \},
\end{align}
where can make the last step by noting that $b(0, k) = 0$ and $\phi(0) = 0$. As we have that $b(n, 1) = \phi(n)$, and that $\calZ_n\{ \phi(n) \} =1- (1-z)^{\alpha} $, we find that
\begin{align}
 \calZ_n\{b(n, k) \} 
&=(1-(1-z)^{\alpha} )^k, \nonumber \\
&= \sum_{\ell=0}^k  \binom{k}{\ell} (-1)^{\ell} (1-z)^{\ell \alpha}, \nonumber \\
&= \sum_{n=0}^\infty (-1)^n \sum_{\ell=0}^k (-1)^{\ell} \binom{k}{\ell}  \binom{\ell\alpha}{n} z^n, \nonumber 
\end{align}
thus we have 
\begin{equation} \label{eq:bnk}
 \bbP_\gamma[ T_k(\cdot) = n ] = b(n,k) = (-1)^n \sum_{\ell = 0}^k (-1)^\ell \binom{k}{\ell} \binom{\ell \alpha}{n} \, .
\end{equation}

We can take similar steps to find a similar relationship as in \eqref{eq:bnk_induction} for $c(n,k)$, that is, we have $\calZ_n\{c(n, k) \}  = \calZ_n\{b(n, k-1)\} \calZ_n\{ \Phi(n) \}$. As $\calZ_n\{ \Phi(n) \} = (1-z)^{\alpha-1}$, we have that 
\begin{align} \label{eq:cnk} 
 \bbP_\gamma[ T_{k-1}(\cdot) < n &< T_{k}(\cdot) ] = c(n,k) \nonumber \\
&= (-1)^n \sum_{\ell = 0}^{k-1} (-1)^\ell \binom{k-1}{\ell} \left [ \binom{(\ell+1)\alpha - 1}{n} - \binom{\ell \alpha}{n}\right] \, .
\end{align}
and thus
\begin{align} \label{eq:P_k} 
 \bbP_\gamma[ T_{k}(\cdot) \le n < T_{k+1}(\cdot) ] &= b(n,k) + c(n,k+1) \nonumber \\
&= (-1)^n 
\sum_{\ell = 0}^{k} (-1)^\ell \binom{k}{\ell}  \binom{(\ell+1)\alpha - 1}{n} 
\, .
\end{align}

\begin{figure}
\centering
\includegraphics[width=\textwidth]{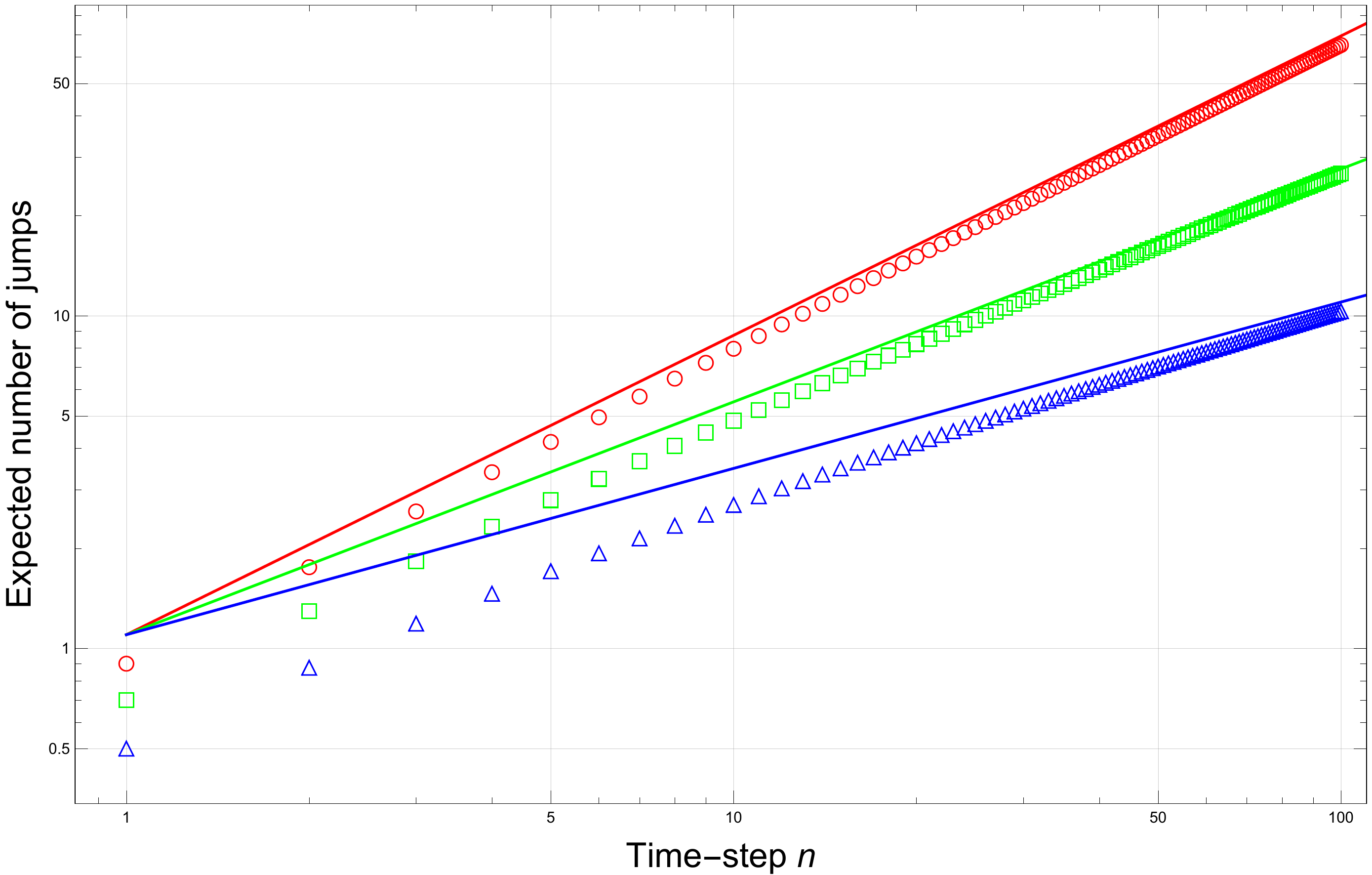} 
\caption{\label{fig:expected_jumps} Expected number of jumps against time step $n$ for $\alpha=0.9$ (red circles), $\alpha=0.7$ (green squares), and $\alpha=0.5$ (blue triangles). The solid lines are plots of $1.1n^\alpha$ to show the asymptotic behaviour.}
\end{figure}

This allows us to calculate the expected value of the number of jumps by the time step $n$,
\begin{align} \label{eq:expected_k}
 \bbE_\gamma[ k : T_k \le n < T_{k+1} ]
&= \sum_{k=0}^\infty k \bbP_\gamma[ T_{k}(\cdot) \le n < T_{k+1}(\cdot) ] \nonumber \\
&= \sum_{k=0}^\infty k \, (-1)^n \sum_{\ell = 0}^{k} (-1)^\ell \binom{k}{\ell}  \binom{(\ell+1)\alpha - 1}{n} \, .
\end{align}
Computationally it is straightforward to  ascertain that 
\begin{equation} \label{eq:exp_bound}
\bbE_\gamma[ k : T_k \le n < T_{k+1} ] = \calO( n^\alpha ).
\end{equation}
This implies that for smaller $\alpha$, we can expect fewer jumps by time step $n$. We plot the expected number of jumps, equation \eqref{eq:expected_k} and the bound in \eqref{eq:exp_bound} in Figure \ref{fig:expected_jumps}.

There is a clear advantage to a path-wise Monte Carlo approach to simulating our discrete subdiffusive process, compared with using the time-stepping approach, solving the mean-field solution $U(i,n)$ using \eqref{eq:sibuya_GME}  \cite{ADHN2015}.
We refer to the the Monte Carlo approach as MC-DTRW and the finite difference approach as
FD-DTRW. The key advantage in the MC-DTRW is the smaller growth of computational complexity as a function of $n$, particularly for small $\alpha$. This is due to the fact that for small $\alpha$ the expected number of jumps up to time step $n$ is small.

Given a space-grid size $L$, we can see from the summation in the left hand side of \eqref{eq:sibuya_GME} that $\calO(L n^2)$ calculations are required to calculate $U(i,n)$ given some initial condition $U(i, 0)$. For some fixed grid spacing $\Delta x, \Delta t$  we have shown in Section \ref{sec:power_law} that $U(\lfloor x / \Delta x\rfloor, \lfloor t / \Delta t\rfloor)$ approximates $u(x,t)$, the solution to \eqref{eq:frac_diffusion}. Thus for a fixed time $t$ and grid spacing $\Delta x, \Delta t$, we know that $\calO(t^2 \Delta t^{-2} \Delta x^{-1})$ computations are required to approximate $u(x,t)$. Now, if we consider $\Delta x$, $\alpha$ and $D_\alpha$ fixed then we must have a time spacing $\Delta t$ determined by the relationship \eqref{eq:delta_t_dep}, that is $\Delta t^\alpha = r \Delta x^2 / (2 D_\alpha)$, and we thus have that the computational complexity of approximating $u(x,t)$ is $\calO(t^2 \Delta x^{-(1 + 4/\alpha)} )$.

By contrast calculating a single path on the same grid requires at most $\calO(n) = \calO(\Delta t^{-1})$ points. In fact, as calculated in \eqref{eq:expected_k} in Section \ref{sec:expected_jumps}, the expected number of jumps to the $n^\mathrm{th}$ time step is $\calO(n^\alpha)$. Hence if we have a fixed $t$, $\Delta x$ and $D_\alpha$, the relationship \eqref{eq:delta_t_dep} gives us that the complexity of calculating a single path is $\calO(n^\alpha) = \calO(t^\alpha \Delta t^{-\alpha}) = \calO(t^\alpha \Delta x^{-2})$ in expectation.

Typically a Monte Carlo procedure would involve $N$ independent trials, or paths. Demonstrating the number of trials necessary for $\varepsilon$ convergences of $U(\lfloor x / \Delta x\rfloor, \lfloor t / \Delta t\rfloor)$ to $u(x,t)$ under some metric is beyond the scope of this paper.

Thus the key advantage of the MC-DTRW method
is that the computational complexity of the method scales as $t^\alpha$ rather than quadratically with $t$ as is the case in the FD-DTRW method.


\section{Discrete time subordination}

The results of the previous section enable us to express the solution of the difference equation \eqref{eq:sibuya_GME} in a closed form. We do so by decomposing the probability of a particle being at grid point $i$, at time $n$, into the probability of being at point $i$ after $k$ jumps, multiplied by the likelihood of making exactly $k$ jumps up to time $n$, for all $k$. We use equation \eqref{eq:P_k} in Section \ref{sec:expected_jumps} for the probability of making $k$ jumps by time $n$. We further note that, for a simple random walk that starts at the origin on an infinite domain, after $k$ jumps the probability of finding the particle at site $i$ is related to the binomial distribution, 
\[
\bbP_\gamma [X_k(\cdot) = i] =\left\{ \begin{array}{c c}\binom{k}{\frac{k+i}{2}} p_r^{\frac{k+i}{2}} p_\ell^{\frac{k-i}{2}}& \mathrm{for}\;\;\frac{k+i}{2}\in \mathbb{Z}^{+}\\ 0 & \mathrm{otherwise}  \end{array}\right. .
\]
Thus we have

\begin{align}
 U(i,n) 
&= \bbP_\gamma [ X(n,\cdot) = i ] \nonumber \\
&= \sum_{k=0}^n \bbP_\gamma [X_k(\cdot) = i] \, \bbP_\gamma[ T_{k}(\cdot) \le n < T_{k+1}(\cdot) ] \nonumber \\
&=\left\{ \begin{array}{c c} \sum_{j=0}^{\frac{n}{2}} \binom{2j}{j+\frac{i}{2}} p_r^{j+\frac{i}{2}} p_\ell^{j-\frac{i}{2}} \,(-1)^n \sum_{\ell = 0}^{2j} (-1)^\ell \binom{2j}{\ell}  \binom{(\ell+1)\alpha - 1}{n}& \mathrm{for}\;i\;\mathrm{even}\\ 
\sum_{j=0}^{\frac{n-1}{2}} \binom{2j+1}{j+\frac{1}{2}+\frac{i}{2}} p_r^{j+\frac{1}{2}+\frac{i}{2}} p_\ell^{j+\frac{1}{2}-\frac{i}{2}} \,(-1)^n \sum_{\ell = 0}^{2j+1} (-1)^\ell \binom{2j+1}{\ell}  \binom{(\ell+1)\alpha - 1}{n}& \mathrm{for}\;i\;\mathrm{odd} \end{array} \right. .
\end{align}
This result bears resemblance to subordination solutions for the continuous time subdiffusion equation \cite{B2001}.

\section{Numerical results} \label{sec:num}

To verify this algorithm we compared the results of the MC-DTRW method outlined in \S \ref{sec:MC}  to the analytic solution of the fractional diffusion equation \eqref{eq:frac_diffusion}, as well as the solution of the discrete master equation \eqref{eq:sibuya_GME} (the FD-DTRW method). 

The results were computed on the domain $x \in [-1,1]$, with zero-flux boundary conditions at $x=1,-1$, i.e. $\frac{\partial u}{\partial x}(-1, t) = \frac{\partial u}{\partial x}(1, t) = 0$. The initial condition was set to the delta function at the origin, $u(x,0) = \delta_0(x)$, equivalent to letting all paths start at the origin for the MC-DTRW simulation, and setting $U(x,0) = \delta_{i,0}$ for the FD-DTRW simulation. 

Using separation of variables it is possible to write a series solution for \eqref{eq:frac_diffusion} subject to these boundary and initial conditions, 
\begin{equation}
u(x,t)=\frac{1}{2}+\sum_{n=1}^{\infty}(-1)^n E_{\alpha}\left(-(n\pi)^2D_\alpha t^{\alpha}\right)\cos(n\pi(x-1))
\end{equation}
The infinite series was approximated by taking the first 900 terms in the sum. This was sufficient to ensure that the results were with to within $10^{-6}$ of the exact solution for the worst case parameters.

\begin{figure}
\centering
\includegraphics[width=1\textwidth]{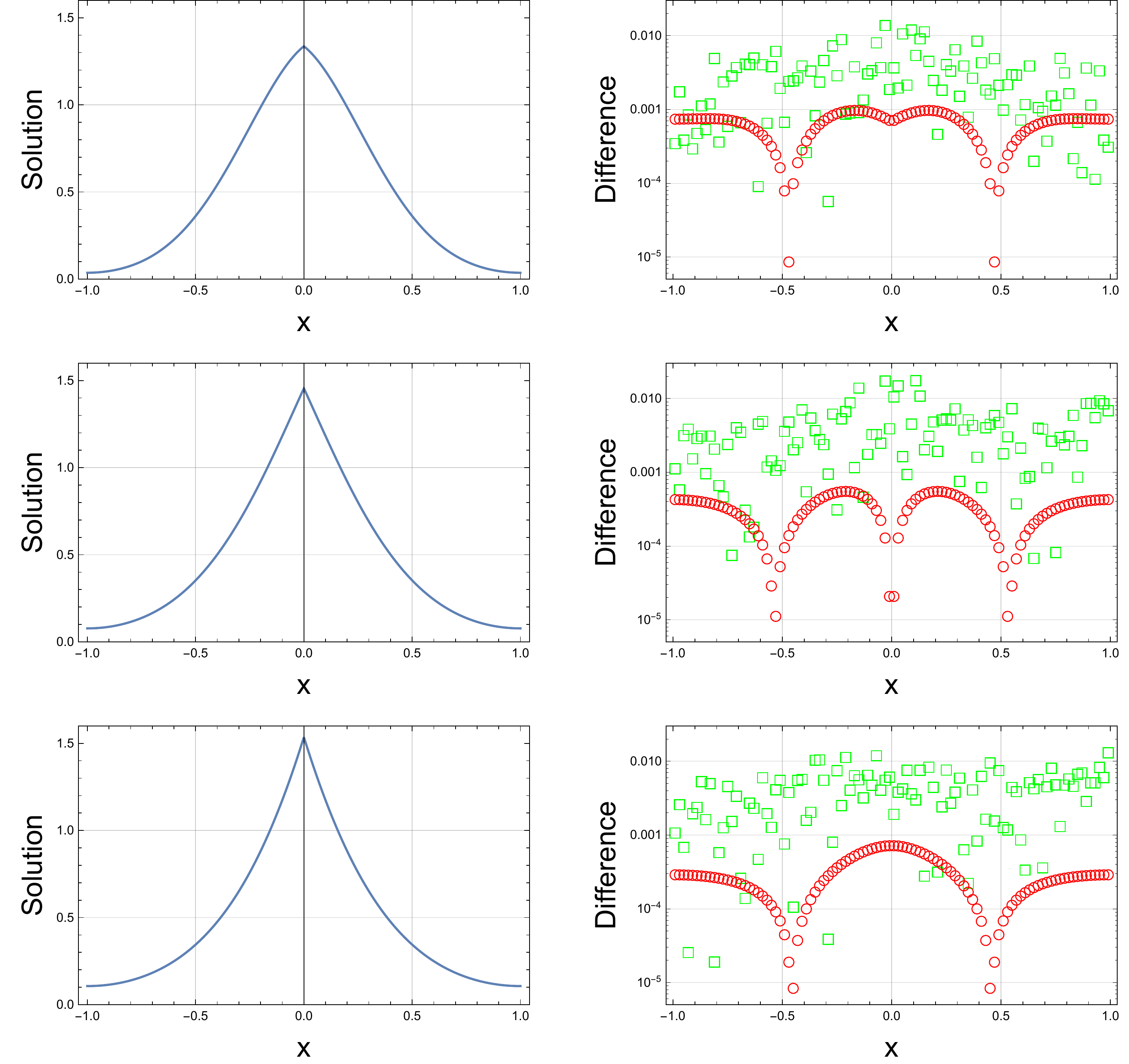}
\caption{ \label{fig:results} Solutions for (top to bottom) $\alpha = 0.9$, $0.7$, and $0.5$, with $N = 10^6$. The right hand column shows the
absolute difference between the analytic solution and MC-DTRW solution (squares) and between the analytic solution and the FD-DTRW solution (circles).}
\end{figure}

We calculated the solutions at time $t=0.5$, with $\alpha = 0.5, 0.6, 0.7, 0.8$ and $0.9$, and $D_\alpha = 0.1$. The grid spacing was set at $\Delta x = 0.2$, with $\Delta t$ determined by equation \eqref{eq:delta_t_dep}. We computed $10^7$ trials for the MC-DTRW simulation. Figure \ref{fig:results} shows the results of the simulations for $\alpha=0.9$, $0.7$ and $0.5$, as well as the absolute difference between the numerical results, both FD-DTRW and MC-DTRW, and the analytic results. 

We see that the difference between the MC-DTRW results and the analytic results was of the order of $10^{-3}$, while the FD-DTRW results were of order $10^{-4}$ from the analytic results. Note that as $N\to\infty$, we would have that the MC-DTRW simulation would converge to the FD-DTRW solution, not the analytic solution.

\begin{figure}
\centering
\includegraphics[width=\textwidth]{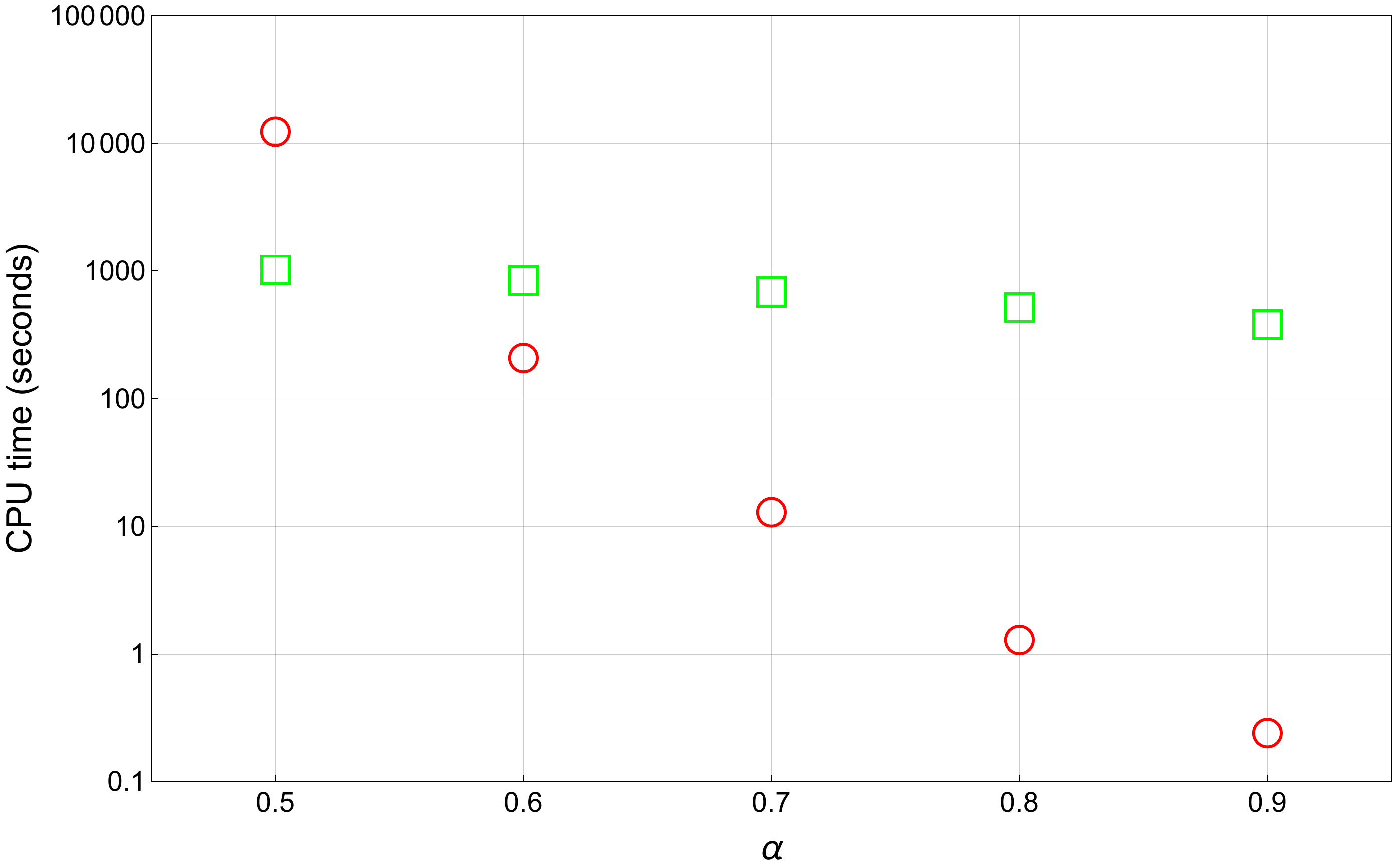}
\caption{ \label{fig:alpha_timing} CPU time comparison against $\alpha$ for the MC-DTRW algorithm (squares) and the FD-DTRW algorithm (circles), $N=10^7$}
\end{figure}

Figure \ref{fig:alpha_timing} shows us the timing of the MC-DTRW simulation compared to the FD-DTRW algorithm for the different $\alpha$. We see that for small $\alpha$, the FD-DTRW takes exponentially more time, while the CPU time for the MC-DTRW solution remains flat. This is as predicted in Section 4. 

\begin{figure} 
\centering
\includegraphics[width=\textwidth]{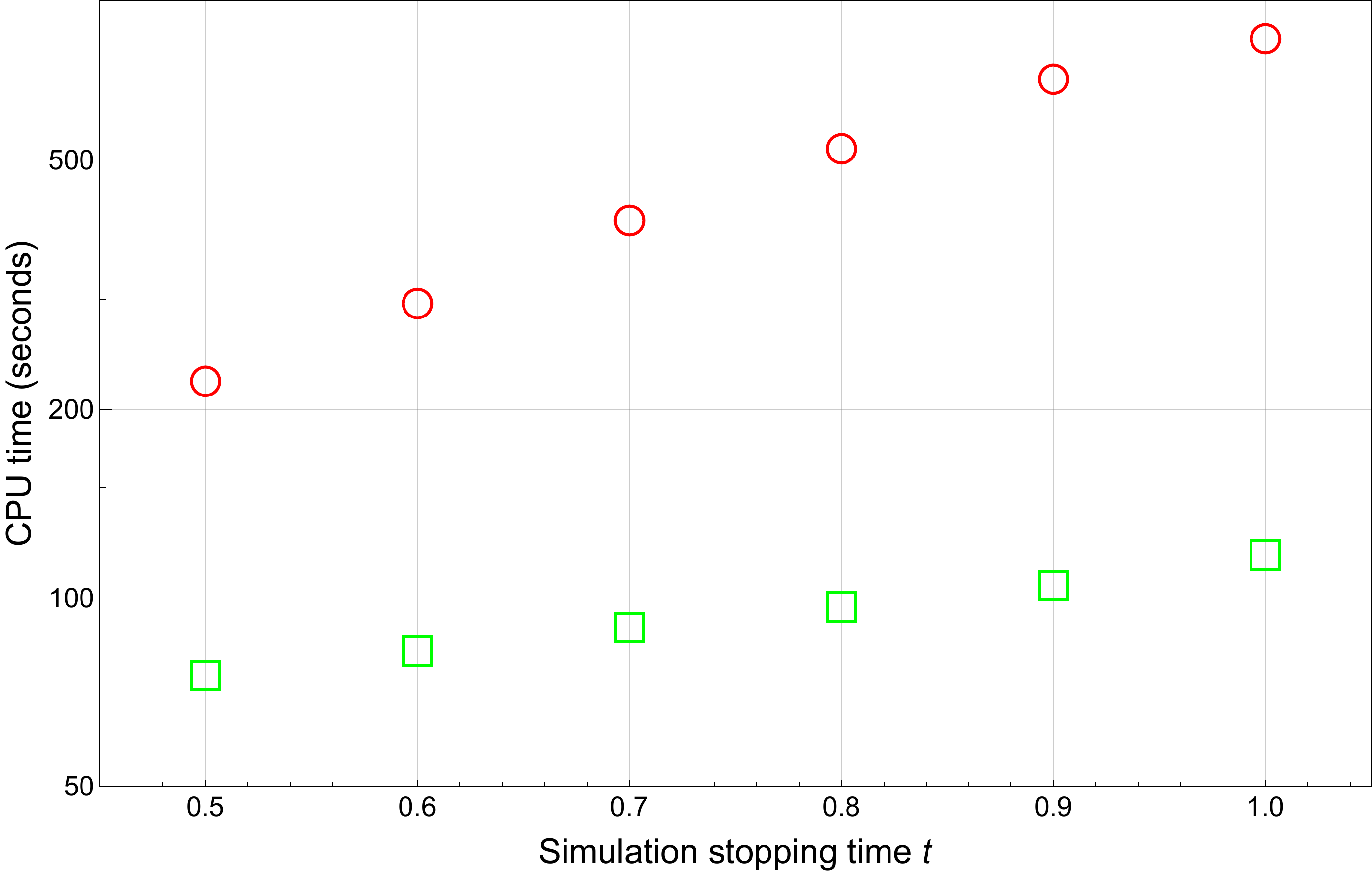}
\caption{\label{fig:t_timing}CPU time comparison against simulation time $t$ for the MC-DTRW algorithm (squares) and the FD-DTRW algorithm (circles), $\alpha=0.7$, $N=10^6$}
\end{figure}

In Figure \ref{fig:t_timing} we examine the CPU time the FD-DTRW and MC-DTRW algorithms took to complete, for various simulation times $t$. Here $\alpha$ was fixed to $0.7$. Here we see the growth in computational complexity against $t$, and in particular how the growth of complexity of the FD-DTRW scheme far outstrips the MC-DTRW simulation.

We note that the CPU time benchmark that we use is highly sensitive to the details of the implementation of our algorithm, and could vary wildly with a slightly different approach. However, we see from our results that the CPU times scale as we would expect from our calculations in Section 4. As the scaling can be so extreme in the case of very small $\alpha$ and large $t$, we find it plausible that no matter what the details of an implementation, an MC-DTRW approach would be wise to consider.

\section{Acknowledgements}
This work was supported by the Australian Commonwealth Government (ARC DP140101193).

\bibliographystyle{elsarticle-num}
\bibliography{dtrw_pathwise}

\end{document}